\documentclass[12pt]{article} 
\usepackage{amsmath, amssymb}
\usepackage{amsthm} 
\usepackage[all]{xy}

\theoremstyle{plain} 

\theoremstyle{definition} 

\def\L{{\cal L}}

\def\bgn{\begin}

\def\L{{\cal L}}
\def\Ob{\text{\rm Ob}}
\def\1{{[1]}}
\def\2{{[2]}}
\def\3{{[3]}}
\def\({\left(}
\def\){\right)}
\def\s-circ{\,{\scriptstyle{\circ}}\,}
\def\<<{<\negthinspace \negthinspace<}
\def\Ad{\text{\rm Ad}}

\def\bgn{\begin}

\def\endaln{\end{align}}

\def\<{<\negthinspace \negthinspace <}

\def\t{\theta}

\def\({\left(}
\def\){\right)}

\def\[{\big[\neg\big[}
\def\]{\big]\neg\big]}
\def\al{\al}

\def\M{{\cal M}}

\def\a{\alpha}
\def\b{\beta}

\def\e{\varepsilon}
\def\gam{\gamma}
\def\Gam{\Gamma}

\def\del{\delta}
\def\lam{\lambda}
\def\Lam{\Lambda}
\def\ome{\omega}
\def\Ome{\Omega}
\def\sig{\sigma}

\def\R{\Bbb R}
\def\C{\Bbb C}
\def\H{\text{\rm H}}
\def\Z{\Bbb Z}
\def\O{\Bbb O}
\def\P{\frak P}
\def\M{\frak M}

\def\w{\wedge}
\def\({\left(}
\def\){\right)}

\def\neg{\negthinspace}
\def\h{\hat}

\def\til{\tilde}

\def\ol{\overline}
\def\pa{\partial}

\def\arrow{\longrightarrow}

\def\bsh{\backslash}

\def\:{\, :\,}

\def\complex{generalized complex }
\def\K\"ahler{generalized K\"ahler}

\def\10{\displaystyle L^{10}}
\def\2{\displaystyle L^2}
\def\c0{\displaystyle C^0}

\def\10{\displaystyle L^{10}}
\def\2{\displaystyle L^2}
\def\del{\delta}
\def\del2{\displaystyle L^2_{0,\delta}}
\def\c0{\displaystyle C^0}

\def\del{\delta}

\def\K{{\cal K}}
\def\M-A{\text{\rm Monge-Amp\`ere}}
\def\O{{\cal O}}

\def\M-A{\text{\rm Monge-Amp\`ere}}
\def\[{\big[\,}
\def\]{\,\big]}
\def\End{\text{\rm End\,}}

\def\P{\Bbb P}
\title{On the stability of locally conformal K\"ahler structures} 
\author{
\textsc{Ryushi Goto$^{*}$} 
}
\date{} 
\begin{document}
\maketitle
\footnote{ 
2010 \textit{Mathematics Subject Classification}.
Primary 53C55; Secondary 32G05.
}
\footnote{ 
\textit{Key words and phrases}. locally conformally K\"ahler metrics, non-K\"ahler surfaces, deformation theory
}
\footnote{ 
$^{*}$Partly supported by the Grant-in-Aid for Scientific Research (C),
Japan Society for the Promotion of Science. 
}
\begin{abstract}
In this article we develop a new approach to the problem of the stability of locally conformally K\"ahler structures (l.c.k structures) under small deformations of complex structures and deformations of flat line bundles. 
We show that under the certain cohomological condition the stability of l.c.k structures does hold. 
We apply our approach to Hopf manifolds 
and its generalizations to obtain the stability of l.c.k structures which do not have potential in general. We give an explicit description of the cohomological obstructions of the stability of l.c.k structures on Inoue surfaces with $b_2=0$.
\end{abstract}
\tableofcontents
\numberwithin{equation}{section}
\section{Introduction}
Let $X$ be a compact complex manifold with complex structure $J$ and $g$ a Hermitian metric on $X$. We denote by $\ome$ the fundamental $2$-form of 
$g$. 
If there is a $d$-closed $1$-form $\eta$ such that $d\ome=\eta\w\ome$, 
then $\ome$ is called {\it a locally conformally K\"ahler structure} (l.c.k structure) on $X$ with {\it Lee form} $\eta$. 
An l.c.k structure yields a K\"ahler metric on the universal covering  of $X$.
Many interesting l.c.k structures on complex manifolds were already constructed.
Hopf surfaces admit l.c.k structures \cite{GO}, \cite{Bel}.
Vaisman metrics \cite{Vai} are l.c.k structures with parallel Lee form, which are constructed
on the quotients of the cones of Sasakian manifolds \cite{BG}.
 Ornea and Vebitsky \cite{OV1} studied the class of l.c.k structures with potential which is the one admitting global K\"ahler potential on the universal covering of $X$.
 Tricerri \cite{Tri} gave l.c.k structures on certain Inoue surfaces with $b_2=0$ \cite{Ino}, which do not have potential. 
 Fujiki and Pontecorvo \cite{FP} used the Twistor theory to provide l.c.k structures on the certain complex surfaces of type VII including hyperbolic Inoue and parabolic Inoue surfaces. 
Brunella constructed l.c.k structures on Kato surfaces which are complex surfaces of type VII with global spherical cell \cite{Bru1},\cite{Bru2}.
\par
Kodaira and Spencer \cite{K.S III} showed that K\"ahler structures are stable under small deformations of complex structures. More precisely, if $X$ is K\"ahlerian, then any small deformation $X_t$ of $X=X_0$
is also K\"ahlerian. 
Contrast to K\"ahler structures,  
l.c.k structures by Tricerri on certain Inoue surfaces with $b_2=0$ are not stable under small deformations of complex structures \cite{Bel}.
However l.c.k structures with potential are stable \cite{OV1}.
This suggests 
the need for further research on the stability of l.c.k structures under deformations.
\par
The purpose of this paper is to obtain the cohomological criterion for the stability of l.c.k structures.
We apply the method developed on deformations of generalized K\"ahler geometry
\cite{Go-1}, \cite{Go0}, \cite{Go1}, \cite{Go2}, \cite{Go3}.
An l.c.k structure $\ome$ gives a flat line bundle $L$ over $X$
which we call {\it the corresponding flat line bundle }to $\ome$.
Then using flat line bundle $L$-valued differential forms $\w^\bullet\otimes L$, 
we have the $L$-valued de Rham cohomology groups $H^\bullet(X, L)$ of 
$L$-valued de Rham complex $(\w^\bullet\otimes L, d_L)$ and the $L$-valued Dolbeault cohomology groups $H^{p,q}(X,L)$ of the 
$L$-valued Dolbeault complex. 
Then it turns out that $\ome$ gives a $d_L$-closed $L$-valued form $\til\ome$
(see Section 1 for more detail).
The $\pa\ol\pa$-lemma for $L$-valued forms is the following:\par\clearpage
{\indent\sc Definition 2.1.}
  Let $X$ be a complex manifold and $L$ a flat line bundle over $X$. 
  Then we say that 
  $(X, L)$ satisfies the $\pa\ol\pa$-lemma at degree $(p,q)$ if 
  there is an $L$-valued form $\gam$ of type $(p, q-1)$ such that 
  $$\pa_L \ol\pa_L \gam=\pa_L\a,
  $$
  for every $\ol\pa_L$-closed $L$-valued form $\a$ of type $(p,q)$.
  \par\medskip
  Then we have the following criterion for the stability of l.c.k structures:\par\medskip
  {\indent\sc Theorem 2.2.}
  Let $X=(M, J)$ be a compact, complex manifold with a locally conformally 
  K\"ahler structure $\ome$ \text{\rm(l.c.k structure)}. 
  We denote by $L$ the corresponding flat line bundle to the l.c.k structure $\ome$.
  We assume that 
  $(X, L)$ satisfies the $\pa\ol\pa$-lemma at degree $(0,2)$.
  Let $\{J_t\}$ be deformations of complex structure $J$ which analytically 
  depend on $t$ and $J_0=J$, where $|t|<\e,$ for a constant $\e>0$. 
 Then there is a positive constant $\e'<\e$ and an analytic family of $2$-forms $\{\ome_t\}$ which satisfies the followings:
 \renewcommand{\labelenumi}{\text{\rm (\roman{enumi})}}
 \bgn{enumerate}
 \item $\ome_0=\ome$,
 \item $\ome_t$ is an l.c.k structure on $(M, J_t)$ for all $t$,
 \item The line bundle $L$ is the corresponding flat line bundle to $\ome_t$
 for all $t$,
 \end{enumerate}
 where $|t|<\e'$.
 \par\medskip
We also consider deformations of flat line bundle $\{L_s\}$ and 
try to construct a family of l.c.k structures $\{\ome_s\}$ such that $L_s$ is the corresponding line bundle to $\ome_s$. 
Then an obstruction to deformations appears as a cohomology class in $H^3(X, L)$
(Theorem \ref{th:obstruction H3}).
We further obtain the criterion for the stability under deformations of flat line bundles:\par\medskip
{\indent\sc Theorem 2.3.}
 Let $\{L_s\}$ be deformations of flat line bundles which analytically depend on $s$ and $L_0=L$, where 
$|s|<\e$ for a constant $\e>0$.
If $(X, L)$ satisfies the $\pa\ol\pa$-lemma at degree $(0,2)$
and $H^3(X, L)=\{0\}$, then there is a positive constant $\e'<\e$ and an analytic family of $2$-forms $\{\ome_s\}$ which satisfies the followings:
\renewcommand{\labelenumi}{\text{\rm (\roman{enumi})}}
 \bgn{enumerate}
 \item $\ome_0=\ome$,
 \item $\ome_s$ is an l.c.k structure on $(M, J)$ for all $s$,
 \item The line bundle $L_s$ is the corresponding flat line bundle to $\ome_s$
 for all $s$,
 \end{enumerate}
where $|s|<\e'$.
\par\medskip
In section 2, we give preliminary results of l.c.k structures and show that our main theorems.  
In section 3, we give proofs of the main theorems. 
In section 4, 5, we discuss the stability of  several examples of l.c.k structures on certain Hopf manifolds and its generalizations. 
They satisfy the $\pa\ol\pa$-lemma at degree $(0,2)$ for a class of flat line bundles.
The $L$-valued Bott-Chern cohomology groups are calculated on them. When the $L$-valued Bott-Chern cohomology group is not trivial,  there are l.c.k structures which do not have potential. For the class of l.c.k structures which do not have potential,
our theorems can be applied to 
obtain deformations of the l.c.k structures. 
In section 6 we give a result on the class of complex surfaces with effective anti-canonical line bundle.
In section 7, we discuss obstructions to deformations of l.c.k structures on Inoue surfaces with $b_2=0$.
There are three classes of Inoue surfaces: $S_M$, $S_{N,p,q,r;t}^{(+)}$ and 
$S_{N,p,q,r}^{(-)}$. It is known that both Inoue surfaces $S_M$ and $S_{N,p,q,r}^{(-)}$ are rigid and 
the Inoue surface $S_{N,p,q,r;t}^{(+)}$ admits $1$-dimensional deformations of complex structures which are parameterized by complex numbers $t \in \C$.
Belgun \cite{Bel} showed that  
 both $S_M$ and $S_{N,p,q,r}^{(-)}$ have l.c.k structures and yet
 $S_{N,p,q,r;t}^{(+)}$ admits an l.c.k structure for only real numbers $t \in \R$.
Thus the stability theorem of l.c.k structures does not hold on $S_{N,p,q,r;t}^{(+)}$.
We show that the obstruction to the stability appears as a non-trivial cohomology class
(Proposition \ref{prop:7.2}).
Further together with Belgun's result, we show that the corresponding flat line bundle to
every l.c.k structure on $S_{N,p,q,r;t}^{(+)}$ must be the canonical line bundle (Proposition \ref{prop:7.3}).
We also show there are obstructions to the stability of l.c.k structures under deformations of flat line bundles on Inoue surfaces with $b_2=0$. 
It is shown that the l.c.k structures on $S_M$, $S_{N,p,q,r;t}^{(+)}$ and $S_{N,p,q,r}^{(-)}$ are not stable under deformations of flat line bundles (Proposition \ref{prop:7.5}, \ref{prop:7.6}).

\section{Stability theorem of locally conformally K\"ahler structures}
Let $M$ be a compact differential manifold of dimension $2n$ with an integrable complex structure $J$. 
We denote by $X$ be the compact complex manifold $(M,J)$.
Let $\ome$ be a locally conformally K\"ahler structure (l.c.k structure) on $X=(M,J)$ with Lee form $\eta$.
We note that if the Lee form $\eta=df$ is $d$-exact, then the metric $e^{-f}g$ is K\"ahlerian. 
We take a covering $\{ U_i\}$ of $M$ such that 
$H^p(U_{i_1}\cap \cdots U_{i_q}, \R)=0$, for all $p>0$ and any 
$i_1, \cdots, i_q $. 
Then $\eta|_{U_i}=df_i,$ where $f_i$ is a function on $U_i$ and 
$f_i-f_j =\lam_{ij}$ is a constant on $U_i \cap U_j$. 
Then $\{e^{-\lam_{ij}}\}_{i \in \Lam}$ gives a class $e^{-\lam}$ of the first 
$\R^*$-valued C\v{e}ck cohomology group $\check{H}^1(M, \R^*)$. 
The class $e^\lam$ yields the flat line bundle $L$ which 
 has a trivialization on $U_i$ with
 locally constant transition functions 
$\{(e^{-\lam_{ij}}, U_i \cap U_j)\}$.
The flat line bundle $L$ over $X=(M,J)$ is called {\it the corresponding flat line bundle to an l.c.k structure $\ome$.}
The vector bundle of $L$-valued $k$-forms is denoted by 
$\w^k \otimes L$. 
An $L$-valued $k$-form $\til\a$ is a section of $\w^k \otimes L$ which is given as a set of $k$-form $\{\a_i\}$ such that 
$\til\a_i=e^{-\lam_{ij}}\til\a_j$ for all $i,j$.
The exterior derivative $d$ induces the differential operator $d_L$ acting on $L$-valued forms 
which yields the $L$-valued de Rham complex: 
$$
\cdots \overset{d_L}{\arrow} \w^k \otimes L \overset{d_L}{\arrow} \w^{k+1}\otimes L \overset{d_L}{\arrow}\cdots.
$$
We denote by $H^\bullet(X, L)$ the $L$-valued de Rham cohomology groups.
Let $d_\eta$ be the differential operator $d-\eta$ which acts on differential forms.
The $\eta$-twisted cohomology groups $H^\bullet_\eta(X)$ are the cohomology groups of the $\eta$-twisted complex:
$$
\cdots \overset{d_\eta}{\arrow} \w^k \overset{d_\eta}{\arrow} \w^{k+1}\overset{d_\eta}{\arrow}\cdots.
$$
A $k$-form $\a$ is regarded as an $L$-valued $k$-form $\til\a:=\{e^{-f_i}\a\}$. 
This gives an isomorphism between complexes $(\w^\bullet, d_\eta)\cong 
(\w^\bullet\otimes L, d_L)$ since $d_\eta=e^{f_i}\circ d \circ e^{-f_i}$ on $U_i$.
Thus we have an isomorphism $H^k_\eta(X)\cong H^k(X, L)$. 
We also denote by $\w^{p,q}\otimes L$ the sheaf of $L$-valued form of type $(p,q)$ on a complex manifold $X=(M, J)$.
Then we have the decomposition $d_L=\pa_L+\ol\pa_L$, where
\bgn{align*}
&\pa_L: \w^{p,q}\otimes L \arrow \w^{p+1, q}\otimes L,\\
&\ol\pa_L:\w^{p,q}\otimes L \arrow \w^{p, q+1}\otimes L.
\end{align*} 
Then we have the Dolbeault complex of $L$-valued forms 
$$
\cdots\arrow \w^{p,q}\otimes L \arrow \w^{p, q+1}\otimes L \arrow \cdots.
$$
We denote by  $H^{p,q}(X, L)$ the $L$-valued Dolbeault cohomology groups. 
The Lee form $\eta$ is decomposed into $\eta^{1,0}$ and $\eta^{0,1}$,
where $\eta^{1,0}\in \w^{1,0}$ and 
$\eta^{0,1}\in \w^{0,1}$.
Let $\pa_\eta$ be the operator $\pa+\eta^{1,0}$ and $\ol\pa_\eta=\ol\pa+\eta^{0,1}$ the its
complex conjugate.
Then we have the twisted Dolbeault cohomology groups $H^{p,q}_\eta(X)$
which are the cohomology groups of the twisted Dolbeault complex: 
$$
\cdots\overset{\ol\pa_\eta}{\arrow}\w^{p,q}\overset{\ol\pa_\eta}{\arrow}\w^{p,q+1}\overset{\ol\pa_\eta}{\arrow}\cdots.
$$
As in before, we have an isomorphism $H^{p,q}_\eta(X)\cong H^{p,q}(X,L)$.
\par
For an l.c.k structure $\ome$, we have a local K\"ahler form $\ome_i=e^{-f_i}\ome$ on $U_i$ and 
a set of K\"ahler form $\til\ome:=\{ \ome_i\}$
gives an $L$-valued $d_L$-closed, positive-definite form $\til\ome$ of type $(1,1)$ on $(M,J)$. 
We call $\til\ome$ an $L$-valued K\"ahler form on $(M.J)$.
Conversely, an $L$-valued K\"ahler form $\til\ome$ yields an l.c.k structure 
$\ome$ and $L$ is the corresponding flat line bundle to $\ome$.

The inclusion $\R^*\to \O_X^*$ induces the map  
$H^1(X, \R^*)\arrow H^1(X, \O_X^*)$. 
Thus the class $e^{-\lam}\in H^1(X, \R)$ gives the holomorphic line bundle $\L$ over $X$. 
Note that $\L$ is a topologically trivial complex line bundle. 
Let $H^p(X, \Ome^q \otimes \L)$ denote the cohomology group, where $\Ome^q$ is the vector bundle of holomorphic $q$ forms. 
Then we have $H^q(X, \Ome^p \otimes \L)\cong H^{p,q}(X,L)$.
\bgn{definition}
  Let $X=(M, J)$ be a complex manifold and $L$ a flat line bundle over $X$. 
  Then we say that 
  $(X, L)$ satisfies the $\pa\ol\pa$-lemma at degree $(p,q)$ if 
  there is an $L$-valued form $\gam$ of type $(p, q-1)$ such that 
  $$\pa_L \ol\pa_L \gam=\pa_L\a,
  $$
  for every $\ol\pa_L$-closed, $L$-valued form $\a$ of type $(p,q)$.
  \end{definition}
  Note that if $H^{p,q}(X,L)=\{0\}$, then $(X,L)$ satisfies the $\pa\ol\pa$-lemma at $(p,q)$.

  \bgn{theorem}\label{th:stability of l.c.k C}
  Let $X=(M, J)$ be a compact, complex manifold with a locally conformally 
  K\"ahler structure $\ome$ \text{\rm(l.c.k structure)}. 
  We denote by $L$ the corresponding flat line bundle to the l.c.k structure $\ome$.
  We assume that 
  $(X, L)$ satisfies the $\pa\ol\pa$-lemma at degree $(0,2)$.
  Let $\{J_t\}$ be deformations of complex structure $J$ which analytically 
  depend on  $t$ and $J_0=J$, where $|t|<\e,$ for a constant $\e>0$. 
 Then there is a positive constant $\e'<\e$ and an analytic family of $2$-forms $\{\ome_t\}$ which satisfies the followings:
 \renewcommand{\labelenumi}{\text{\rm (\roman{enumi})}}
 \bgn{enumerate}
 \item $\ome_0=\ome$,
 \item $\ome_t$ is an l.c.k structure on $(M, J_t)$ for all $t$,
 \item The line bundle $L$ is the corresponding flat line bundle to $\ome_t$
 for all $t$,
 \end{enumerate}
 where $|t|<\e'$.
  \end{theorem}
Next we consider deformations of the flat line bundle $L$ on a complex manifold $X=(M,J)$ with an l.c.k structure $\ome$. 
 Small deformations of the flat line bundle $L$ are parametrized by the first cohomology group $H^1(X, \R)$. 
 Let $\{L_s\}$ be deformations of flat line bundles with $L_0=L$ which are parametrized by the real number $s$. 
 Then deformations $\{L_s\}$ is given by a family of $d$-closed $1$-forms 
 $\{\eta_s\}$ which are representatives of $H^1(M,\R)$.
We have the following criterion for the stability of l.c.k structures under deformations of flat line bundles
 \bgn{theorem}\label{th:stability of l.c.k L}
Let $\{L_s\}$ be deformations of 
flat line bundles which analytically depend on $s$ and $L_0=L$, where 
$|s|<\e$ for a constant $\e>0$.
If $(X, L)$ satisfies the $\pa\ol\pa$-lemma at degree $(0,2)$
and 
$H^3(X, L)=\{0\}$,  then 
there is a positive constant $\e'<\e$ and an analytic family of $2$-forms 
$\{\ome_s\}$ which satisfies the followings:
\renewcommand{\labelenumi}{\text{\rm (\roman{enumi})}}
\bgn{enumerate}
\item $\ome_0=\ome$,
\item $\ome_s$ is an l.c.k structure on $(M, J)$ for all $s$,
\item The line bundle $L_s$ is the corresponding flat line bundle to $\ome_s$
for all $s$,
 \end{enumerate}
where $|s|<\e'$.
\end{theorem}
For a family of $d$-closed $2$-forms $\eta_s$, 
we denote by $\dot{\eta}_0:=\frac{d}{ds}{\eta_s}|_{s=0}$ the infinitesimal tangent of deformations of 
 $\{L_s\}$ at $s=0$ which gives the class $[\dot{\eta_0}]\in H^1(X, \R)$. 
 Then the class $[\dot{\eta_0}\w\ome]\in H^3(X, L)$ is regarded as the first obstruction 
 to the existence of a smooth family of l.c.k structures $\{\ome_s\}$ such that
 $L_s$ is the corresponding flat line bundle to  $\ome_s$.
 \bgn{theorem}\label{th:obstruction H3}
 If the class $[\dot{\eta_0}\w\ome]\in H^3(X, L)$ does not vanish, then 
 $X$ does not admit a smooth family of l.c.k forms 
 $\{\ome_s\}$ such that $L_s$ is the corresponding flat line bundle to  $\ome_s$ for all small $s$.
 \end{theorem}
 In Section 7, we shall show that the obstruction does not vanish for
 Inoue surfaces with $b_2=0$. 
 Combining the above two theorems, we obtain the following theorem:
 \bgn{theorem}\label{th:stability of l.c.k CL}
 Let $\{J_t\}$ be deformations of complex structure $J$ which analytically 
 depend on  $t$ and $J_0=J$, where $|t|<\e,$ for a constant $\e>0$.
 We denote by $\{L_s\}$ deformations of 
  flat line bundles which analytically depend on $s$ and $L_0=L$, where 
  $|s|<\e$.
  If $(X, L)$ satisfies the $\pa\ol\pa$-lemma at degree $(0,2)$
and $H^3(X, L)=\{0\}$,
then 
there is a positive constant $\e'<\e$ and an analytic $2$-parameter family of $2$-forms 
  $\{\ome_{s,t}\}$ which satisfies the followings:
  \renewcommand{\labelenumi}{\text{\rm (\roman{enumi})}}
 \bgn{enumerate}
 \item $\ome_{0,0}=\ome$,
 \item $\ome_{s,t}$ is an l.c.k structure on $(M, J_t)$ for all $t$,
 \item The line bundle $L_s$ is the corresponding flat line bundle to 
 $\ome_{s,t}$
 for all $s$,
 \end{enumerate}
where $|s|,\, |t|<\e'$.
\end{theorem}

\section{Proof of main theorems}
This section is devoted to proof theorems \ref{th:stability of l.c.k C}, 
\ref{th:stability of l.c.k L}, \ref{th:obstruction H3} and \ref{th:stability of l.c.k CL}.
 We use the same notation as in the previous section.
 We have already explained that an l.c.k structure $\ome$ on $X=(M,J)$ gives an $L$-valued K\"ahler form $\til\ome$ on $X=(M,J)$ and conversely, an $L$-valued K\"ahler form $\ome$ yields an l.c.k structure 
 $\til\ome$ which admits the corresponding flat line bundle $L$. 
Thus Theorem \ref{th:stability of l.c.k C} 
is reduced to the following theorem:
\bgn{theorem}\label{th:stability of L.K C}
  Let $X=(M, J)$ be a compact, complex manifold with an $L$-valued 
  K\"ahler structure $\til\ome$. 
  We assume that $(X, L)$ satisfies the $\pa\ol\pa$-lemma at degree $(0,2)$.
  Let $\{J_t\}$ be deformations of complex structure $J$ which analytically 
  depend on  $t$ and $J_0=J$, where $|t|<\e,$ for a constant $\e>0$. 
 Then there is a positive constant $\e'<\e$ and an analytic family of $L$-valued $2$-forms $\{\til\ome_t\}$ which satisfies the followings:
 \renewcommand{\labelenumi}{\text{\rm (\roman{enumi})}}
 \bgn{enumerate}
 \item $\til\ome_0=\til\ome$,
 \item $\til\ome_t$ is an $L$-valued K\"ahler form on $(M, J_t)$ for all $t$,
 \end{enumerate}
 where $|t|<\e'$.
  \end{theorem}
From now on, we identify a locally conformally K\"ahler structure with the corresponding $L$-valued K\"ahler form and 
for simplicity, we denote by $\ome$ an $L$-valued K\"ahler form on $X=(M,J)$, in stead of $\til\ome$ in this section.
\par
The bundle of complex endmorphism End$(TM)$ is decomposed into two components:
$E_0:=(\w^{1,0}\otimes T^{0,1})\oplus(\w^{0,1}\otimes T^{1,0})$ and 
$E_1:=(\w^{1,0}\otimes T^{1,0})\oplus (\w^{0,1}\otimes T^{0,1})$. 
We denote by $E_i^\R$ the real part of $E_i$ for $i=0,1.$ 
An endmorphism $a$ acts on the complex structure $J$ by 
the adjoint $ad_aJ:= [a, J]$ at each point.
The exponential $e^a$ is regarded as an element of GL$(TM)$ for $a \in $End$(TM)$ which acts on 
the complex structure $J=J_0$ by the Adjoint action
$\Ad_{e^a}J$. 
The exponential $e^a$ also acts on the $L$-valued K\"ahler form $\ome$ by  $e^a \cdot\ome$ which is 
the linear representation of the Lie group GL$(TM)$. 
We also denote by $a \cdot\ome$ the action of End$(TM)$ on $\ome$ which is the action of the Lie algebra.
We assume that there is a family of complex structures $\{J_t\}$ as in Theorem \ref{th:stability of L.K C}. 
Then there is a family of endmorphisms 
$\{a_t\}$ with $a_t\in E^\R_0$ such that 
\bgn{equation}
J_t=\Ad_{e^{a_t}}J,
\end{equation} 
where $a(t)$ is an analytic family (cf. Proposition 2.6 page 541 in \cite{Go1}), 
$$
a(t) =a_1t+a_2\frac{t^2}{2!}+\cdots.
$$
Since $\Ad_{e^b}J=J$ for $b \in E_1^\R$, the action of $E_1^\R$ preserves the complex structure $J$.
We assume that there is a family of endmorphisms $\{b_t\}$ with $b_t \in E^\R_1$ such that 
\bgn{equation}\label{eq:dL=0}
d_L(e^{a_t}\circ e^{b_t}\cdot \ome)=0
\end{equation}
By the Campbel-Hausdorff formula, there is a section $z(t)\in $End$(TM)$ such that 
$$
e^{z(t)}=e^{a(t)}\circ e^{b(t)}.
$$
Then it turns out that 
$e^{z(t)}\cdot\ome$ is an $L$-valued K\"ahler form on $(M, J_t)$ for each $t$. 
In fact, the diagonal action of $e^{z(t)}$ on the pair $(J, \ome)$ yields a family of pairs 
$(\Ad_{e^{z(t)}}J,\,\,\, e^{z(t)}\cdot\ome)$, where $e^{z(t)}\cdot\ome $ is an $L$-valued Hermitian form of type $(1,1)$ with respect to the complex structure $\Ad_{e^{z(t)}}J$.
Then we have 
\bgn{align}
\Ad_{e^{z(t)}}J =&\Ad_{e^{a(t)}}\circ \Ad_{e^{b(t)}}J \\
=&\Ad_{e^{a(t)}}J =J_t.
\end{align}
Thus it follows from (\ref{eq:dL=0}) that $(J_t, e^{z(t)}\cdot\ome)$ is an $L$-valued K\"ahler form on $(M, J_t)$.
We shall construct a formal power series $\{b_t\}$ with $b_t\in E_1^\R$ which satisfies (\ref{eq:dL=0}).
The action of $E^\R_1$ on $\ome$ gives $L$-valued real forms of type $(1,1)$ and the action of $E^\R_0$ on $\ome$ yields $L$-valued real forms of type $(0,2)$ and type $(2,0)$. 
Since $\ome$ is a non-degenerate $2$-form at each point, it turns out that 
every form of type $(1,1)$ is written as $b \cdot\ome$ for a section $b \in E^\R_1$. Thus  
$L$-valued real forms $\w^{1,1}_\R\otimes L$ of type $(1,1)$ are given by 
$\{\, b \cdot\ome \, |\, b \in E^\R_1\, \}$ and
we also have 
$(\w^{2,0}\oplus\w^{0,2})_\R \otimes L =\{a \cdot\ome\, |\, a \in E_0^R\}$, 
where $(\w^{2,0}\oplus\w^{0,2})_\R$ is the real part of $\w^{2,0}\oplus\w^{0,2}$. 
\bgn{lemma}\label{key lemma}
We assume that $(X, L)$ satisfies the $\pa\ol\pa$-lemma at degree $(0,2)$. 
If $d_L\a$ is a real $d_L$-exact form of $(\w^{2,1}\oplus\w^{1,2})_\R\otimes L$ for a real $\a\in\w^2\otimes L$, 
then there is a real $L$-valued form $\b$ of type $(1,1)$ such that 
$d_L \b= d_L\a$, where $\b$ is written as 
$b \cdot\ome$ for $b \in E_1^\R$.

\end{lemma}
\bgn{proof}We decomposed $\a$ by 
$\a =\a^{2,0}+\a^{1,1}+\a^{0,2}$, where $\a^{p,q}\in \w^{p,q}\otimes L$ and 
$\a^{2,0}=\ol{\a^{0,2}}$, since $\a$ is real.
Since $d_L\a$ is of type $(2,1)$ and $(1,2)$, we have $\ol\pa_L\a^{0,2}=0$ and $\pa_L\a^{2,0}=0$.
Since $(X, L)$ satisfies the $\pa\ol\pa$-lemma, 
there is an $L$-valued form $\gam$ of type $(0, 1)$ such that 
  $$\pa_L \ol\pa_L \gam=\pa_L\a^{0,2}.
  $$
  Then it follows that $\ol\pa_L \pa_L \ol\gam=\ol\pa_L\a^{2,0}.$
  We define $\b \in \w^{1,1}_\R\otimes L$ by 
 $\b=-\pa_L \gam-\ol\pa_L \ol{\gam}+\a^{1,1}$. 
 Then $\b$ satisfies that 
$$
d_L\b=d_L\a.
$$
\end{proof}
In order to obtain an estimate of $\b$ which satisfies $d_L\b=d_L\a$, 
we use the Hodge theory of the following$L$-valued complex: 
\bgn{equation}\label{eq:complex}
 \w^{1,1}_\R\otimes L \overset{d_L}\arrow 
(\w^{2,1}\oplus\w^{1,2})_\R\otimes L \overset{d_L}\arrow 
(\w^{3,1}\oplus\w^{2,2}\oplus\w^{1,3})_\R\otimes L \overset{d_L}\arrow \cdots.
\end{equation}
The complex (\ref{eq:complex}) is elliptic at the terms $(\w^{2,1}\oplus\w^{1,2})_\R\otimes L$
and 
$(\w^{3,1}\oplus\w^{2,2}\oplus\w^{1,3})_\R\otimes L$.
We denote by $d_L^*$ be the formal adjoint of $d_L$ and by 
$G$ the Green operator which acts on the term $(\w^{2,1}\oplus\w^{1,2})_\R\otimes L$, 
where we take a metric.
Then we define $\b \in \w^{1,1}_\R\otimes L$ by 
$$
\b=d_L^*Gd_L\a.
$$
Then it follows form Lemma \ref{key lemma} that $d_L\b =d_L\a$ and by using the elliptic estimate, we have
$$
\|\b\|_s \leq C\|\a\|_{s},
$$
where $\|\, \|_s$ denotes the Sobolev norm for $s>0$.

\indent{\sc Proof of Theorem \ref{th:stability of L.K C}.}
We shall obtain a power series $b(t)$ by the induction on the degree $k$ of $t$.
The first term $a_1 \in E^\R_0$ is decomposed into ${a_1}'$ and $a_1''$, where 
$a_1' \in \w^{0,1}\otimes T^{1,0}$ and $a_1''=\ol{a_1'}$.
The component $a_1'$ is $\ol\pa$-closed which gives the Kodaira-Spencer class $[a_1']\in H^1(X, T^{1,0})$ of deformations $\{J_t\}$. 
Let $d_L(e^{z(t)}\cdot\ome)_{[k]}$ denotes the term of $d_L(e^{z(t)}\cdot\ome)$ of degree $k$ in $t$.
The first term is given by 
\bgn{equation}\label{eq:the first obs}
d_L(e^{z(t)}\cdot\ome)_{[1]}:= d_L(a_1 \cdot\ome)+ d_L(b_1 \cdot\ome)=0
\end{equation}
Since $\ol\pa a_1'=0$ and $a_1''=\ol{a_1'}$, we have 
$d_L(a_1 \cdot\ome) \in (\w^{2,1}\oplus\w^{1,2})_\R\otimes L$. 
Then it follows from lemma \ref{key lemma} that 
there is a $b_1 \in E^\R_1$ such that $d_L(a_1 \cdot\ome)+ d_L(b_1 \cdot\ome)=0$.
Next we consider an operator $e^{-z(t)}\circ d_L \circ e^{z(t)}$ acting on $L$-valued differential forms. 
It follows that the operator $e^{-z(t)}\circ d_L \circ e^{z(t)}$ is locally written as a composition of the Lie derivative and the interior,  the exterior product 
\bgn{equation}\label{eq:Clifford-Lie operator}
e^{-z(t)}\circ d_L \circ e^{z(t)}=\sum_j {\t^j}\w\L_{v_J} +N_j,
\end{equation}
where $\t^j$ is the exterior product by a $1$-form $\t^j$ and 
$\L_{v_j}$ is the Lie derivative by the vector $v_j$ and $N_j$ is 
an element of $\w^2T^*M \otimes TM$ which acts on differential forms by the interior and  the exterior product
(cf. \cite{Go-1} Lemma 2.4 and 2.7, page 221-223 and Definition 2.2 in \cite{Go0} for more detail).
We take that an open covering $\{ U_\a\}$ of $M$ and a nowhere-vanishing holomorphic $n$-form $\Ome_\a$ on each $U_\a$. We denote by $\Phi_\a$ the pair $(\Ome_\a, \ome)$. 
Then the Lie derivative of $\Phi_\a$ by a vector $v_j$ is given by 
$$
\L_{v_j}\Phi_\a= (a_\a \cdot\Ome_\a, \,\,\, a_\a \cdot\ome),
$$
for $a_\a \in \End(TM)=\w^1 \otimes TM$.
Thus it follows from (\ref{eq:Clifford-Lie operator}) that there is a section $h_\a \in \w^2T^*M \otimes TM$ such that 
$$
e^{-z(t)}\circ d_L \circ e^{z(t)}\cdot \Phi_\a=
(h_\a \cdot\Ome_\a, \,\,h_\a \cdot\ome).
$$
Since $J_t=\Ad_{e^{z(t)}}J$, 
$e^{z(t)}\cdot\Ome_\a$ is a form of type $(n,0)$ with respect to the complex structure $J_t$.
 Since $J_t$ is integrable, $d \,e^{z(t)}\cdot\Ome_\a$ is a form of type $(n,1)$ with respect to the complex structure $J_t$.
Thus we have that $e^{-z(t)}\circ d_L \circ e^{z(t)}\cdot\Ome_\a \in\w^{n,1}$.
It follows that the component of $h_\a$ in  $\w^{0,2}\otimes T^{1,0}$ must vanish. 
Thus the component of $h_\a \cdot\ome$ of $\w^{3,0}$ also vanish and it implies that 
$e^{-z(t)}\circ d_L \circ e^{z(t)}\cdot\ome=h_\a \cdot\ome$ is of type $(2,1)$ and $(1,2)$. 
We assume that we already have $b_1, \cdots, b_{k-1}\in E^\R_1$ such that 
$d_L  \(e^{b(t)}\cdot\ome \)_{[i]}=0$, 
for all $0 \leq i\leq k-1$, where 
$b(t) =\sum_i b_i \frac{t^i}{i!}$. 
Since $e^{z(t)}=e^{a(t)}\circ e^{b(t)}$, we have 
\bgn{align*}
d_L (e^{z(t)}\cdot\ome )_{[k]} =&\sum_{\stackrel {i+j=k}{i,j \geq 0}}
\(e^{-z(t)} \)_{[i]} d_L \(e^{z(t)}\cdot \ome \)_{[j]} \\
=&\(e^{-z(t)}\circ d_L \circ e^{z(t)} \cdot\ome \)_{[k]}
\end{align*}
Then it follows that $d_L (e^{z(t)}\cdot\ome )_{[k]}=(h_\a \cdot\ome)_{[k]}$ is of type $(2,1)$ and $(1,2)$. 
The $k$th term $d_L (e^{z(t)}\cdot\ome)_{[k]} $ is written as the sum of the linear term and the non-linear term $\Ob_k$
$$
d_L (e^{z(t)}\cdot\ome )_{[k]} =\frac 1{k!} d_L( b_k \cdot\ome ) +\Ob_k,
$$
where $\Ob_k$ depends only on $a(t)$ and $b_1, \cdots, b_{k-1}$ which is $d_L$-exact.
Thus it follows from lemma \ref{key lemma} that there is a $b_k \in E^\R_1$ such that 
$d_L (e^{z(t)}\cdot\ome )_{[k]} =\frac1{k!} d_L( b_k \cdot\ome) +\Ob_k=0$.
By our assumption of the induction, we have a solution $b(t)$ in the form of the formal power series.
It turns out that the $b(t)$ is a convergent series which gives a smooth family 
of $L$-valued K\"ahler forms $\{\ome_t\}=\{e^{z(t)}\cdot\ome \}$ by the same method as in  \cite{Go-1, Go1, Go2, Go3}.
\qed
\\
Theorem \ref{th:stability of l.c.k L}
is equivalent to the following:

\bgn{theorem}\label{th:stability of L-l.c.k L}
 Let $X=(M, J)$ be a compact, complex manifold with an $L$-valued 
  K\"ahler structure $\ome$. 
 Let $\{L_s\}$ be deformations of 
  flat line bundles which analytically depend on $s$ and $L_0=L$, where 
  $|s|<\e$ for a constant $\e>0$.
  If $(X, L)$ satisfies the $\pa\ol\pa$-lemma at degree $(0,2)$ and 
  $H^3(X, L)=\{0\}$,  then 
  there is a positive constant $\e'<\e$ and an analytic family of $L_s$-valued $2$-forms 
  $\{\ome_s\}$ which satisfies the followings:
  \renewcommand{\labelenumi}{\text{\rm (\roman{enumi})}}
 \bgn{enumerate}
 \item $\ome_0=\ome$,
 \item $\ome_s$ is an $L_s$-valued K\"ahler on $(M, J)$ for all $s$,
 \end{enumerate}
where $|s|<\e'$.
\end{theorem}

\indent{\sc Proof of theorem \ref{th:stability of L-l.c.k L} and Theorem \ref{th:obstruction H3}. }
Deformations of flat line bundles $\{L_s\}$ are given by a family of elements in $H^1(X, \R)$ which are written as a family of smooth $d$-closed $1$-forms $\{\eta_s\}$ with $\eta_0=0$.
Then $d_L +\eta_s$ is the differential operator and a flat section of $L_s$ is given by 
a section $\sig$ of $L$ with 
 $(d_L+\eta_s)\sig=0$.
As before, we see that an $L_s$-valued K\"ahler form is given by 
an $L$-valued Hermitian form $\ome_s$ such that $(d_L+\eta_s)\ome_s=0$.

Thus we shall obtain a family of section $b(s)$ with $b(s)\in E^\R_1$ such that 
$(d_L+\eta_s)(e^{b(s)}\cdot\ome)=0$. 
The first term of the equation is given by 
$$
d_L b_1 \cdot\ome + \eta_1 \cdot\ome =0.
$$
Thus the class $[\eta_1\cdot\ome]=[\eta_1\w\ome]\in H^3(X, L)$ must vanish if there is a family of 
$L_s$-valued K\"ahler forms. 
Hence Theorem \ref{th:obstruction H3} is proved.

From our assumption in Theorem \ref{th:stability of L-l.c.k L}, we have  $H^3(X, L)=\{0\}$. 
Thus it follows that 
$\eta_1 \w\ome$ is $d_L$-exact real form of type $(2,1)$ and $(1,2)$. 
Then applying Lemma \ref{key lemma}, we obtain $b_1 \in E^\R_1$ such that 
$d_L b_1 \cdot\ome + \eta_1 \cdot\ome =0.$
We shall use the induction on $k$. 
We assume that we already have $b_1,\cdots, b_{k-1}$ such that 
$(d_L+\eta_s)(e^{b(s)}\cdot\ome)_{[i]} =0$ for all $0 \leq i <k$. 
Then we have 
$$
(d_L+\eta_s)(e^{b(s)}\cdot\ome)_{[k]}=\(e^{-b(s)}\circ (d_L+\eta_s)\circ e^{b(s)}\cdot\ome \)_{[k]} 
$$
Since $b(s)\in E^\R_1$, there is a $h_\a$ on each open set $U_\a$ as before such that
$e^{-b(s)}\circ d_L \circ e^{b(s)}\cdot\ome =h_\a \cdot\ome\in 
(\w^{2,1}\oplus \w^{1,2})_\R\otimes L$. 
Further $e^{-b(s)}\circ \eta_s \circ e^{b(s)}$ is also a $1$-form and then 
we see that 
$(d_L+\eta_s)(e^{b(s)}\cdot\ome)_{[k]}\in (\w^{2,1}\oplus \w^{1,2})\otimes L$. 
The term $(d_L+\eta_s)(e^{b(s)}\cdot\ome)_{[k]}$ is given by  
$$(d_L+\eta_s)(e^{b(s)}\cdot\ome)_{[k]}=\frac1{k!}d_L b_k \cdot\ome +\Ob_k(\eta_s).$$
It follows from our assumption that 
$\(d_Le^{b(s)}\cdot\ome\)_{[j]}=-\(\eta_s\w e^{b(s)}\cdot\ome\)_{[j]}$ for 
$j=0,1,\cdots, k-1$.
Since $\eta_s$ is a $d$-closed $1$-form and its degree is greater than or equal to one, we have
\bgn{align}
d_L\circ (d_L+\eta_s)\(e^{b(s)}\cdot\ome\)_{[k]}=&
-\(\eta_s\w d_Le^{b(s)}\cdot\ome\)_{[k]}\\
=&-\sum_{\stackrel{i+j=k}{i\geq 1}}\(\eta_s\)_{[i]}\,\(d_Le^{b(s)}\cdot\ome\)_{[j]}\\
=&\sum_{\stackrel{i+j=k}{i\geq 1}}\(\eta_s\)_{[i]}\w\(\eta_s\w e^{b(s)}\cdot\ome\)_{[j]}\\
=&\sum_{i+k}\(\eta_s\w\eta_s\w e^{b(s)}\)_{[k]}=0.
\end{align} 
Thus $\Ob_k(\eta_s)$ is $d_L$-closed which gives a class of $H^3(X, L)$. 
From our assumption $H^3(X, L)=\{0\}$, $\Ob_k(\eta_s)$ is a $d_L$-exact form of type $(2,1)$ and $(1,2)$. 
Then applying Lemma \ref{key lemma}, we have $b_k \in E^\R_1$ such that 
$\frac1{k!}d_L b_k +\Ob_k(\eta_s)=0$. 
Thus from our assumption of the induction, we have a solution $b(s)$ in the form of the formal power series which can be shown to be a convergent series. We obtain a smooth family of $L_s$-valued K\"ahler forms $\{\ome_s\}$.
\qed

Theorem \ref{th:stability of l.c.k CL} is also equivalent to the following:
\bgn{theorem}\label{th:stability of L-l.c.k CL}
 Let $X=(M, J)$ be a compact, complex manifold with an $L$-valued 
  K\"ahler structure $\ome$. 
 Let $\{J_t\}$ be deformations of complex structure $J$ which analytically 
  depend on  $t$ and $J_0=J$, where $|t|<\e,$ for a constant $\e>0$.
 We denote by $\{L_s\}$ deformations of 
  flat line bundles which analytically depend on $s$ with $L_0=L$, where 
  $|s|<\e$.
  If $(X, L)$ satisfies the $\pa\ol\pa$-lemma at degree $(0,2)$
  and 
  $H^3(X, L)=\{0\}$,  then 
  there is a positive constant $\e'<1$ and an analytic $2$-parameter family of 
  $L_s$-valued $2$ form 
  $\{\ome_{s,t}\}$ which satisfies the followings:
  \renewcommand{\labelenumi}{\text{\rm (\roman{enumi})}}
 \bgn{enumerate}
 \item $\ome_{0,0}=\ome$,
 \item $\ome_{s,t}$ is an $L_s$-valued K\"ahler on $(M, J_t)$ for all $s$,
 \end{enumerate}
where $|s|,\, |t|<\e'$.
\end{theorem}

\indent{\sc Proof of Theorem \ref{th:stability of L-l.c.k CL}.}
We shall use the same method as before. 
Deformations of complex structures $\{J_t\}$ is given by a family of endmorphisms
$\{a_t\}$ with $a_t \in E^\R_0$ and deformations of flat line bundles $\{L_s\}$ is also described by 
a family of $d$-closed $1$-forms
$\{\eta_s\}$. 
Then we shall construct a $2$-parameter family of $\{b(s,t)\}$ with $b(s,t)\in E^\R_1$ such that 
$(d_L+\eta_s)( e^{a(t)}\circ e^{b(s,t)}\cdot\ome)=0$. 
The first term of the equation in $t$ is given by 
$$
d_La_1 \cdot\ome + d_Lb_1 \cdot\ome +\eta_1 \w\ome=0
$$
Since $d_L a_1 \cdot\ome \in \w^{2,1}\oplus\w^{1,2}$ and 
$[\eta_1 \w\ome]=0 \in H^3(X, L)$ and $\eta_1 \w\ome\in (\w^{2,1}\oplus\w^{1,2})_\R$, we have a solution $b_1 \in E_1^\R$ of the first equation from Lemma \ref{key lemma} 
We also have a solution $b_k$  of the $k$th term of the equation and the solution $b(t)$ by the same method 
as in the proof of Theorem \ref{th:stability of L.K C}. 
\qed
\bgn{remark}
If an l.c.k structure is given as an $L$-valued K\"ahler form with potential:
$$
\ome=\sqrt{-1}\pa_L \ol\pa_L \phi,
$$
for an $L$-valued function $\phi$, then obstructions to deformations vanish. 
In fact, we can solve the equation (\ref{eq:the first obs}) explicitly.
it follows from $\ol\pa a_1'=0$ and $a_1''=\ol{a_1'}$ that we have 
\bgn{align}
d_L (a_1 \cdot\ome)=&\sqrt{-1}d_L(a_1'+a_1'')\cdot\pa_L \ol\pa_L \phi \\
=&\sqrt{-1}d_L(\pa_La_1' \pa_L \phi-\ol\pa_La_1'' \ol\pa_L \phi).
\end{align}
Thus $b_1=\sqrt{-1}(\pa_La_1' \pa_L \phi-\ol\pa_La_1'' \ol\pa_L \phi)
\in \w^{1,1}_\R\otimes L$ gives a solution of (\ref{eq:the first obs}). 
This is consistent with the result by Ornea and Verbitsky \cite{OV1} that 
l.c.k structures with potential are stable under small deformations of complex structures.
\end{remark}
\bgn{remark}
It is worth to note that our method gives a proof of the stability of K\"ahler structures which is different from the original proof by Kodaira-Spencer.
The method by Kodaira-Spencer depends on the fact that a K\"ahler form is harmonic and 
they applied the harmonic projection to obtain the stability of K\"ahler structures. 
An l.c.k structure is, however not harmonic and we can not apply the method by 
Kodaira-Spencer to obtain the stability of l.c.k structures. 
\end{remark}
\numberwithin{equation}{section}
\section{Hopf manifolds}
\def\P{{\cal P}}
We recall basic facts on Hopf manifolds. A contraction $f: (\C^n,0)\to (\C^n,0)$ is an automorphism of $\C^n$ fixing $0$ with the property that the eigenvalues $\a_1,\cdots, \a_n$ of $f'(0)$, the differential of $f$ at $0$ are inside the unit circle. 
A diagonal  contraction is an automorphism of the form $f: (z_1,\cdots, z_n)\to 
(\a_1z_1,\cdots,\a_n z_n)$, where $0<|\a_1|<\cdots<|\a_n|<1$.

A Hopf manifold is defined to be the quotient of $\C^n-\{0\}$ by the group of the automorphisms generated by a contraction $f$.
A Hopf manifold is a compact complex manifold which is diffeomorphic to $S^1\times S^{2n-1}$.

A diagonal Hopf manifold is the quotient by the group of the automorphisms generated by a diagonal contraction $f$.
A classical Hopf manifold is a diagonal Hopf manifold given by the contraction 
of the type $f:(z_1,\cdots, z_n) \to (\a z_1, \cdots,\a z_n),$
$0<|\a|<1$.
 It was shown that diagonal Hopf manifold admits an l.c.k metric. 
(see further \cite{GO}, \cite{OV1} and \cite{Bel} for l.c.k structures on Hopf manifolds). 
In the special case of  $\a_1=\cdots \a_n=\a \in \R$, an l.c.k metric on the diagonal Hopf $X$ is explicitly written as 
$$
\ome_\lam =\frac{-1}{2i}\frac{\pa\ol\pa \gam^\lam}{\gam^\lam},
$$
where $\lam >0$.
and $\gam^2=\sum_{i=1}^n |z_i|^2$. 
The Lee form $\eta_\lam$ of $\ome_\lam$ is given by 
$$
\eta_\lam =-\lam \frac{d\gam}{\gam}.
$$
We denote by $d_\eta$ the differential operator $d-\eta$. Then 
we have $d_{\eta_\lam}\ome_\lam =0$. 
We also define operators $\pa_{\eta_\lam}$ and $\ol\pa_{\eta_\lam}$ 
respectively
by 
$$
\pa_{\eta_\lam}=\pa-\eta^{1,0}, \quad \ol\pa_{\eta_\lam}=
\ol\pa-\eta^{0,1},
$$
where we are using the decomposition $\eta=\eta^{1,0}+\eta^{0,1}$ and 
$\eta^{1,0}$ is a form of type $(1,0)$ and $\eta^{0,1}=\ol{\eta^{1,0}}$. 
Let $\rho_\lam: \pi_1(X)\to \R^\times$ be the representation of 
$\pi_1(X)\cong \{\a^n\,|\,n \in \Z\}$ which is given by 
$$
\rho_\lam (\a^n)\mapsto \a^{n \lam},
$$
where $\lam \in \Z$.
Then the real flat line bundle $L_\lam$ over $X$ is defined to be the quotient in terms of the representation $\rho_\lam$: 
$$
L_\lam =\C^n \bsh \{0\}\times_{\rho_{\lam}}\R.
$$
We also denote by $\L_\lam$ the holomorphic line bundle given by the flat line bundle $L_\lam$.
Then we see that $\ome_\lam$ gives the $L_\lam$-valued K\"ahler form $\til\ome_\lam$: 
$$
\til\ome_\lam =\frac{-1}{2i}\pa\ol\pa \gam^\lam.
$$
The $L_\lam$-valued Dolbeault cohomology groups on the diagonal Hopf manifold $X$ was calculated in \cite{Ise}, \cite{Mall}. 
In particular, it was shown that $H^{0,2}(X, L_\lam)=\{0\}$ for $\lam>0$. 
Thus we have that every l.c.k structure with $L_{\lam}$ ($\lam>0$) on 
the diagonal Hopf manifold $X$ is stable under small deformations of complex structures of $X$. 
 We note that our result can be applied to l.c.k structures on $X$ which do not have any potential in general.
In fact, the $L_\lam$-valued Bott-Chern cohomology group $H^1(X,\P(L_\lam))$ of $X$ does not vanish for a positive integer $\lam$ on the classical Hopf manifold, 
where the $L_\lam$-valued Bott-Chern cohomology group $H^1(X, \P(L_\lam))$ is given by 
$$
H^1(X, \P(L_\lam))=\frac{\{\, \a \in \w^{1,1}_\R\, |\, d_{\eta_\lam}\a=0\}}
{\{\pa_{\eta_\lam}\ol\pa_{\eta_\lam} f|f \in C^\infty(X, \R)\}}.
$$
Then we have 
$$
\dim H^1(X,\P(L_\lam))=\dim H^0(\Bbb P^{n-1}, \O(\lam))=
\bgn{pmatrix}
\lam+n-1\\ n-1
\end{pmatrix}
$$
for a positive integer $\lam$
(see the next section for the proof).
Thus $\ome_\lam+\a$ for sufficiently small representative $\a \in H^1(\P(L))$ is an l.c.k which does not have any potential. 
We shall generalize these results to l.c.k structures on the certain quotient of (semi-negative) principal fibre spaces over K\"ahler manifolds with fibre elliptic curves in the next section.

\section{Certain generalizations of Hopf manifolds}
\subsection{$L_\lam$-valued Dolbeault cohomology groups}
Let $E$ be a negative complex line bundle over a compact K\"ahler manifold $B$.
We assume that the anti-canonical line bundle $-K_B$ of $B$ is semi-positive, that is, 
the curvature form of a Hermitian connection of the $-K_B$  over a K\"ahler manifold $B$ is semi-positive.
We note that  Fano manifolds and Calabi-Yau manifolds satisfy the assumption.
We define a compact complex manifold $X$ to be the quotient of 
the complement $E \bsh \{0\}$ by the group of automorphisms generated by 
the multiplication of the constant $\a$ with $\a \neq 0,1$ on each fibre.
Thus we have a principal fibre $\pi: X \to B$ with fibre elliptic curve $\C^*/\Z$ , that is, the group $\C^*/\Z$
acts freely on $X$ and
$B$ is the quotient by the action of $\C^*/\Z$.
Let $\sig$ be a local nonzero holomorphic section of the line bundle $E$ and write 
$\gam:=\|\sig \|_h^2$, where $h$ is a Hermitian metric in $E$. 
Then the $1$-st Chern form is given by 
$-i \pa\ol\pa\log\gam$. 
An l.c.k metric $\ome_\lam$ on $X$ is given by 
$$
\ome_\lam= \frac{-1}{2i}\frac{\pa\ol\pa \gam^\lam}{\gam^\lam},
$$
for a positive $\lam$. 
Since $E$ is a negative line bundle, we see that $\ome_\lam$ is positive-definite.
The Lee form of $\eta_\lam$ of the l.c.k form $\ome_\lam$ is given by 
$$
\eta_\lam =-\lam \frac{d \gam}{\gam}.
$$
We denote by $d_\eta$ the differential operator $d-\eta$. Then 
we have $d_{\eta_\lam}\ome_\lam =0$. 
We also define operators $\pa_{\eta_\lam}$ and $\ol\pa_{\eta_\lam}$,
respectively
by 
$$
\pa_{\eta_\lam}=\pa-\eta^{1,0}, \quad \ol\pa_{\eta_\lam}=
\ol\pa-\eta^{0,1},
$$
where we are using the decomposition $\eta=\eta^{1,0}+\eta^{0,1}$ and 
$\eta^{1,0}$ is a form of type $(1,0)$ and $\eta^{0,1}=\ol{\eta^{1,0}}$. 
Let $\rho_\lam: \pi_1(X)\to \R^\times$ be the representation of 
$\pi_1(X)\cong \{\a^n\,|\,n \in \Z\}$ which is given by 
$$
\rho_\lam (\a^n)\mapsto \a^{n \lam},
$$
where $\lam \in \Z$.
Then the real flat line bundle $L_\lam$ over $X$ is defined to be the quotient by the representation $\rho_\lam: \Z \to \R^\times$: 
$$
L_\lam =E \bsh \{0\}\times_{\rho_{\lam}}\R.
$$
Let $\h\rho_\lam:\C^\times\to \C^\times$ be the representation given by $\h\rho_\lam(w)=w^\lam$ and 
$L_{B,\lam}$ the complex line bundle which is defined by 
$$
L_{B,\lam}=E \bsh \{0\}\times_{\h\rho_{\lam}}\C.
$$
We denote by $\L_{B,\lam}$ the sheaf of holomorphic sections of $L_{B,\lam}$.
Then we see that $\L_\lam =\pi^*L_{B,\lam}$ since $\h\rho_\lam$ is 
restricted to $\rho_\lam$.

Then we see that $\ome_\lam$ gives the $L_\lam$-valued K\"ahler form $\til\ome_\lam$: 
$$
\til\ome_\lam =\frac{-1}{2i}\pa\ol\pa \gam^\lam.
$$
\bgn{proposition}\label{prop:5.1}
$$H^{0,k}(X, L_\lam)=H^k(X,\L)=\{0\}, \qquad (k \geq 2)$$
\end{proposition}
\bgn{proof}
In order to obtain cohomology groups $H^k(X, \L)$,
we apply the Leray spectral sequence of the fibre bundle $\pi:X \to B$.
The $E_2$-term is given by 
$$
E_2^{p,q}=H^p(B, R^q \pi_*\L_\lam), \quad(k=p+q).
$$

By using the projection formula, it follows from $\L_\lam =\pi^*L_{B,\lam}$ that 
$$
R^q \pi_*\L_\lam =R^q \pi_*\O_X \otimes \L_{B,\lam}.
$$
Since the fibre is an elliptic curve, $R^q \pi_*\L_\lam=0$ for $q \neq 0, 1$. 
The form $\eta^{0,1}$ on $X$ is $\ol\pa$-closed which gives a global nowhere-vanishing section of $R^1 \pi_*\O_X$. 
It implies that $R^1 \pi_*\O_X=\O_B$.
Thus it follows that the Leray spectral sequence degenerates at the $E_2$-term.
it follows that $L_{B,-\lam}$ is negative for $\lam>0$. 
Since $-K_B$ is semi-positive, 
$K_B \otimes L_{B,-\lam}$ is negative also for $\lam>0$.
Thus applying the Kodaira vanishing theorem, we have that 
$H^k(B, \L_{B,\lam})\cong H^{n-1-k}(B, \L_{B,-\lam}\otimes K_B)=\{0\}$ and 
$H^{k-1}(B,\L_{B,\lam})\cong H^{n-k}(B,\K_{B,-\lam}\otimes K_B)=\{0\}$, for 
$k=2,\cdots, n.$
(Note that $\dim B=n-1$.)
Thus we have $H^k(X,\L_{\lam})=\{0\}$ for $k \geq 2$.
\end{proof}
\subsection{Bott-Chern cohomology groups}
A real function $f$ on $X$ is pluriharmonic if $\pa\ol\pa f=0$.
Let $\cal P$ be the sheaf of real pluriharmonic functions. 
Then we have the short exact sequence:
\bgn{equation}\label{eq:pluri exact seq}
0 \arrow \underline\R\arrow \O_X \arrow {\cal P}\arrow 0,
\end{equation}
where
$\underline \R$ denotes the real constant sheaf which is a subsheaf of $\O_X$ and the surjection 
$\O_X \to \P$ is given by taking the imaginary part of a holomorphic function, i.e., $\phi\mapsto \frac1{2i}(\phi-\ol \phi)$.
Let $L$ be a real flat line bundle over $X$. 
 The tensor product by $L$ of the short exact sequence (\ref{eq:pluri exact seq}) induces 
a short exact sequence:
$$
0 \arrow L \arrow \L \arrow \P \otimes L \arrow 0
$$
We put $\P(L):=\P \otimes L$. Then we have the long exact sequence: 
\bgn{equation}\label{long Bott-Chern}
\cdots\arrow H^1(X, L) \arrow H^1(X, \L)\arrow H^1(X, \P(L))\arrow H^2(X, L)\arrow \cdots
\end{equation}
Let $\w^{p,q}(L)$ be the sheaf of $L$-valued $C^\infty$ forms of
type $(p,q)$ on $X$. We denote by $\w^{p,p}_\R(L)$ the real part of $\w^{p,p}(L)$.
Then we have the exact sequence of sheaves:
$$
0 \arrow \P(L)\arrow \w^{0,0}_\R(L)
\overset{i \pa\ol\pa}\arrow \w^{1,1}_\R(L)\overset{d}\arrow
\(\w^{2,1}(L)\oplus \w^{1,2}(L)\)_\R \overset{d}\arrow\cdots,
$$
where $\(\w^{2,1}(L)\oplus \w^{1,2}(L)\)_\R$ denotes the real part
of $\(\w^{2,1}(L)\oplus \w^{1,2}(L)\)$.
Thus the cohomology group $H^1(X, \P(L))$ coincides with the real Bott-Chern cohomology group:
$$
H^1(X, \P(L))=\frac{\{\, \a \in \Gam(X,\w^{1,1}_\R(L))\, |\, d_L\a=0\}}
{\{i \pa\ol\pa f|f \in \Gam(X, \w^{0,0}_\R(L)\}},
$$
where $\Gam(X, \w^{1,1}_\R(L))$ 
is the set of $L$-valued real $C^\infty$ forms of type $(1,1)$ and 
$\Gam(X, \w^{0,0}_\R(L))$ is the set of $L^\R$-valued $C^\infty$  real functions 
on $X$ $^{*1}$\footnote{ $*1$ The author is grateful to Prof. Fujiki for this description of the Bott-Chern cohomology groups}.

\bgn{definition}
A Vaisman metric is an l.c.k metric $g$ whose Lee form is parallel with respect to the $g$.
\end{definition}
Let $H^\bullet(Y, L)$ be the cohomology groups of the complex 
$(\w^\bullet\otimes L, d_L)$ on a compact differential manifold $Y$ with a real flat line bundle $L$ over $Y$.
It was shown \cite{LLMP} that $H^k(Y, L)$ vanishes for all $k$ on a
locally conformally symplectic manifold $Y$ if $Y$ admits a compatible Riemannian metric on which the Lee form is parallel.  
Thus it is pointed out in \cite{OV1} that 
\bgn{theorem}\cite{LLMP}
Let $Y$ be a compact complex manifold with a Vaisman metric. Then we have 
$$H^k(Y, L)=\{0\}, \qquad \text{\rm for all }k$$
\end{theorem}
Thus from the long exact sequence (\ref{long Bott-Chern}), we obtain 
\bgn{proposition}\label{prop:5.4}
Let $Y$ be a compact complex manifold with a Vaisman metric.
Then the $L$-valued Bott-Chern cohomology group is given by 
$$H^1(Y, \P(L))=H^1(Y, \L)\equiv H^{0,1}_{\ol\pa}(Y, L)$$,
where $H^{0,1}_{\ol\pa}(Y, L)$ denotes the $L$-valued Dolbeault cohomology group. 
\end{proposition}
We shall apply the proposition to the principal fibre bundle $\pi: X \to B$ as in
Subsection 5.1.
We already show that the Leray spectral sequence degenerates at the $E_2$-term and 
$R^1 \pi_*\O_X$ is trivial (see Proposition \ref{prop:5.1}). The $E_2$-term is given by 
$$
E_2^{1,0}\cong \equiv H^1(B,\L_{B,\lam}), \quad
E_2^{0,1}\cong  H^0(B, \L_{B,\lam})
$$
We already see that $\lam>0$, i.e., $L_{B,\lam}$ is a positive line bundle over $B$. 
Then applying the Kodaira vanishing theorem, we have $H^1(B,\L_{B,\lam})
\equiv H^{n-2}(B, K_B \otimes \L_{B,-\lam})=\{0\}$. 
Thus we have $H^1(X, \L)\cong H^0(B,\L_{B,\lam})$. 
Then it follows from Proposition \ref{prop:5.4} that  
\bgn{proposition}
Let $X$ be the complex manifold as before. 
Then we have 
$$
H^1(X, \P(L_\lam))\equiv H^1(X, \L_\lam)\cong  H^0(B, L_{B,\lam}).
$$
\end{proposition}
\bgn{remark}
In our case that $M$ is a fibre bundle $\pi:M \to B$ with fibre $S^1$ and $\eta$ gives a generator of 
the first cohomology group of each fibre, we directly see that $H^k(X,L)=\{0\}$. 
In fact, applying the Leray spectral sequence to $p: M \to B$, we obtain that 
the $E_2$-terms vanishes since $R^k \pi_* L=\{0\}$ for all $k$.
\end{remark}
As before, the cohomology group $H^1(X, \P(L_\lam))$ is written as 
$$
H^1(X, \P(L))=\frac{\{\, \a \in \w^{1,1}_\R\, |\, d_{\eta_\lam}\a=0\}}
{\{\pa_{\eta_\lam}\ol\pa_{\eta_\lam} f|f \in C^\infty(X, \R)\}}.
$$
Thus as in the previous subsection, 
$\ome_\lam+\a$ is an l.c.k form for a sufficiently small representative $\a$
of $H^1(X, \P(L))$. If the class $[\a]$ does not vanish, the l.c.k structure 
$\ome_\lam+\a$ does not have any potential. 
We denote by $J$ the complex structure of $X$. 
Since $H^3(X, L_\lam)=\{0\}$ and $H^2(X, \L_\lam)=\{0\}$,
we have the stability theorem of l.c.k structures on $X$.
\bgn{proposition}
Let $\ome:=\ome_\lam+\a$ be the l.c.k structure on $X$, where $\ome_\lam$
is a Vaisman metric and $\a$ is a representative of the $L_\lam$-valued Bott-Chern class. 
Let $\{L_s\}$ be deformations of 
 of flat line bundles and $\{J_t\}$ deformations of complex structures of $X$
  which analytically depend on parameters $s, t$, respectively and $L_0=L$, $J_0=J$, where $|s|, |t|<\e$, for a constant $\e>0$.
 Then there is a positive constant $\e'$ and an analytic $2$-parameter family of $2$-forms $\{\ome_{s,t}\}$
such that $\ome_{s,t}$ is an l.c.k structure on $(X, J_t)$ with the corresponding flat line bundle $L_s$ for $|s|,\,|t|<\e'$.
\end{proposition}

\section{Complex surfaces with effective anti-canonical line bundle}
Let $\ome$ be  a locally conformally K\"ahler structure on a compact complex surface $S$ and $L$ the corresponding flat line bundle to $\ome$.
We denote by $\L$ the holomorphic line bundle given by $L$.
By the Serre duality, 
we have $H^{2}(S,\L)\cong H^0(S, K \otimes \L^{-1})$. 
\bgn{proposition}
We assume that the anti-canonical line bundle $-K$ is effective and $\L^{-1}$ does not admit any 
holomorphic section. 
Then an l.c.k structure $\ome$ is stable under small deformations of $S$, that is, every small deformation of $S$ admits a locally conformally K\"ahler 
structure.
\end{proposition}
\bgn{proof}
It suffices to show that $H^2(S, \L)=\{0\}$.
Since $-K$ is effective, $K$ is an ideal sheaf and
we have the injective map $H^0(S, K \otimes \L^{-1})\to H^0(S, \L^{-1})$. 
Since $H^0(S, \L^{-1})=\{0\}$, we have $H^0(S, K \otimes \L^{-1})\cong H^2(S, \L)=\{0\}$. 
Thus the result follow from theorem \ref{th:stability of l.c.k L}.
\end{proof}

\section{Inoue surfaces with $b_2=0$}
\def\H{\Bbb H}
These non-K\"ahler surfaces were introduced by Inoue \cite{Ino} which are called Inoue surfaces. 
An Inoue surface $S$ satisfies the following two conditions : 
\renewcommand{\labelenumi}{\text{\rm (\roman{enumi})}}
\bgn{enumerate}
\item
The first Betti number $b_1(S)$ is equal to one and the second Betti number $b_2(S)$ vanishes,
\item $S$ contains no curves.
\end{enumerate}
There are three kinds of Inoue surfaces: $S_M$, $S_{N,p,q,r;t}^{(+)}$ and 
$S_{N,p,q,r}^{(-)}$, all of them being compact quotients of $\Bbb H \times \C$ by 
discrete groups of holomorphic automorphisms, 
where $\Bbb H$ denotes the upper half plane of the complex numbers $\C$. Inoue surface $S_{N,p,q,r;t}^{(+)}$ admits a $1$-dimensional family of deformations which are parameterized by $t \in \C$. 
On the other hand, $S_M$ and $S_{N,p,q,r}^{(-)}$ are rigid.
Tricerri \cite{Tri} constructed locally conformally K\"aher structures on 
$S_M$, $S_{N,p,q,r;t}^{(+)}$ and $S_{N,p,q,r}^{(-)}$ for 
$t \in \R$. 
Belgun \cite{Bel} showed that Inoue surface $S_{N,p,q,r;t}^{(+)}$ does not admit any l.c.k metric for $t \in \C \bsh \R$. 
It implies that l.c.k metrics on Inoue surfaces $S_{N,p,q,r;t}^{(+)}$, $(t \in \R)$ 
are not stable under small deformations. 
We shall show that there is an obstruction to the stability of the l.c.k metrics on $S_{N,p,q,r;t}^{(+)}$. 
\subsection{Inoue surfaces $S_{N,p,q,r;t}^{(+)}$}
Inoue surface $S_{N,p,q,r;t}^{(+)}$ is the quotient 
of $\H \times \C$ by discrete group $G_{N,p,q,r;t}^{(+)}$
$$
S_{N,p,q,r;t}^{(+)}=\H \times \C/ G_{N,p,q,r;t}^{(+)}
$$
We denote by $(w, z)$ holomorphic coordinates of $\H \times \C$, where 
$w=w_1+\sqrt{-1}w_2 \in \H$ and $z=z_1+\sqrt{-1}z_2 \in \C$. 
The group $G_{N,p,q,r;t}^{(+)}$ is generated by the following analytic automorphisms:
\bgn{align}
&\phi_0:(w,z)\mapsto (\a w, z+t)\\
&\phi_i: (w,z)\mapsto (w+a_i, z+b_i w+c_i), \quad i=1,2\\
&\phi_3:(w,z)\mapsto \(w, z+\frac{(b_1a_2-b_2a_1)}{r} \)
\end{align}
where $a_i, b_i , c_i , r$ are real constants which are obtained by choosing a unimodular matrix $N=(n_{ij})\in $SL$(2, \Z)$ with 
two real eigenvalues $\a>1$ and $\frac1\a$,
and two real eigenvectors $(a_1, a_2)$ and $(b_1, b_2)$ corresponding to 
$\a$ and $\frac1\a$, respectively, and  three integers $p,q,r$ ($r \neq 0$).
Two constants $c_1$ and $c_2$ are defined to be a solution of the following equation:
$$
(c_1, c_2) =(c_1, c_2)N^t+(e_1, e_2)+\frac{(b_1a_2-b_2a_1)}{r}(p,q),
$$
where $e_i=\frac12 n_{i,1}(n_{i,1}-1)a_1b_1+\frac12 n_{i,2}(n_{i,2}-1)a_2b_2+n_{i,1}n_{i,2}b_1a_1$, for $i=1,2$.
Note that constants $\a, a_i, b_i, c_i, r$ are real. 
Then the action of $G_{N,p,q,r;t}^{(+)}$ on $\H\times\C$ is properly discontinuous and has no fixed points.
We have a basis $\{\t^1,\, \,\t^2\}$ of invariant forms of type $(1,0)$ under the action of $G_{N,p,q,r;t}^{(+)}$:
\bgn{align}
&\t^1=\frac{dw}{w_2}, \qquad \t^2=\frac{z_2}{w_2} dw-dz,
\end{align}
Note that $\t^1$ and $\t^2$ are not holomorphic.
We also have a basis of invariant vectors of type $(1,0)$ under the action of 
$G_{N,p,q,r;t}^{(+)}$:
\bgn{align}
&X_1=w_2 \frac{\pa}{\pa w}+z_2 \frac{\pa}{\pa z},\qquad 
X_2=-\frac{\pa}{\pa z}
\end{align}
Then we have 
\bgn{align}
d \t^1=&-\frac{dw_2}{w_2}\w \t^1=-\frac{\t^1-\ol\t^1}{2i}\w \t^1=\frac1{2i}\ol\t^1\w\t^1\\
d\t^2=&\t^1 \w\frac1{2i}(\t^2 -\ol\t^2)
\end{align}
We define $\ome$ by 
$\ome=\frac1{i}(\t^1 \w\ol\t^1+\t^2\w\ol\t^2)$ . Then we have 
\bgn{align}
d\ome=&\frac{dw_2}{w_2}\w\ome 
\end{align}
Thus $\ome$ is an l.c.k structure with Lee form $\eta=\frac{dw_2}{w_2}$ \cite{Tri}. 
We denote by $d_\eta$ the operator $d-\eta$. Then $d_\eta\ome =0$.
We define $\til\ome$ by 
$\til \ome =\frac1{w_2}\ome$. 
Then we have 
$\phi_0^*\til\ome=\a^{-1}\til\ome$, 
$\phi_i^*\til\ome=\til\ome$ for $i=1,2,3$ and
$d \til\ome=0$.
Thus $\til\ome$ is a flat line bundle $L$-valued K\"ahler form. 
A holomorphic $2$-form $\Ome:=dw \w dz$ satisfies 
$\phi_0^*\Ome=\a \Ome$ and $\phi_i^*\Ome=\Ome$ for 
$i=1,2,3$. 
Thus $\Ome$ gives a nowhere-vanishing $L^{-1}$-valued holomorphic $2$-form.
Hence we see that $K=L$.
For simplicity, we denote by $S_t^{(+)}$ the Inoue surface $S_{N,p,q,r;t}^{(+)}$.
Then we obtain 
\bgn{lemma}
$H^{0,2}(S_{t}^{(+)},\, L)\cong \C$.
\end{lemma}
\bgn{proof}
Since $K=L$, we have 
$$H^{0,2}(S_{t}^{(+)}, L)\cong H^0(S_{t}^{(+)}, K \otimes L^{-1})\cong H^0(S_{t}^{(+)}, \O)\cong \C.$$
\end{proof}
The Inoue surface $S_{t}^{(+)}$ admits a $1$-dimensional deformations of complex structures and then Kodaira-Spencer class is given by 
$$[X_2\otimes\ol\t^1]=[\frac{\pa}{\pa z}\otimes\frac{d\ol w}{w_2}]\in H^1(S_t^{(+)},\Theta)
$$
We put $\e=X_2\otimes\ol\t^1$. 
The following complex is a key point of deformations of l.c.k structures 
which is equivalent to the complex (\ref{eq:complex}). 
The obstruction to the stability of l.c.k structures arises as a cohomology class 
at the term $(\w^{2,1}\oplus\w^{1,2})_\R$:
\bgn{equation}\label{d eta complex}
\w^{1,1}_\R\overset {d_\eta}{\arrow} (\w^{2,1}\oplus\w^{1,2})_\R\overset {d_\eta}\arrow(\w^{3,1}\oplus\w^{2,2}\oplus\w^{1,3})_\R\overset {d_\eta}{\arrow}\cdots
\end{equation}
The K-S class $[X_2\otimes\ol\t^1]\in H^1(S_t^{(+)}, \Theta)$ acts on $\ome$,
\bgn{align}
(X_2\otimes\ol\t^1)\cdot\ome=
&\frac1{i}(\ol\t^1\w\ol\t^2)
\end{align}
Then we have 
\bgn{align}
d_\eta \((X_2\otimes\ol\t^1)\cdot\ome \) 
=&\frac2{i(w_2)^2}dw_2\w d\ol w\w d\ol z\\
=&\t^1\w\ol\t^1\w\ol \t^2\\
=&i\t^1\w\ol\t^1\w i(dz_1- idz_2)
\end{align}
Then the real part and the imaginary part of $d_\eta(\e\cdot\ome)$
give classes of the cohomology group of the complex (\ref{d eta complex}), respectively. 
\bgn{proposition}\label{prop:7.2}
The class $[i\t^1\w\ol\t^1\w dz_1]$ vanishes, that is,  the form $i\t^1\w\ol\t^1\w  dz_1$ is written as 
$d_\eta\a$ for $\a\in \w^{1,1}_\R$. 
However the class $[i\t^1\w\ol\t^1\w dz_2]$ does not vanish,
where we are considering the cohomology classes of the complex (\ref{d eta complex}).
\end{proposition}
This reflects the result by Belgun that the Inoue surface $S_{N,p,q,r;t}^{(+)}$
admits l.c.k metrics for $t\in \R$, however does not admit l.c.k metric
for $t\in \C\bsh\R$.
\bgn{proof}
Since  we have
$$
d_\eta(\t^1\w\ol\t^2)=d_\eta(\ol\t^1\w\t^2)=
\frac{1}{i}\t^1\w\ol\t^1\w dz_1,
$$
the class $[i\t^1\w\ol\t^1\w  dz_1]$ vanishes.
We apply the Hodge theory to the complex (\ref{d eta complex}). 
Let $(d_\eta)^*$ be the formal adjoint of the differential operator 
$\w^{1,1}_\R\overset{d_\eta}\arrow(\w^{2,1}\oplus\w^{1,2})_\R$ in terms of the Hermitian metric $\ome$. 
Then the formal adjoint $(d_\eta)^*$ is given by 
$$
(d_\eta)^*=-\pi_{\w^{1,1}}\circ (*^{-1} d_{-\eta} *),
$$
where $\pi_{\w^{1,1}}$ denotes the projection to forms of type $(1,1)$ and 
note that $d_{-\eta}=d+\eta$.
Thus we have $(d_\eta)^*\(i\t^1\w\ol\t^1\w  dz_2\)=0$. 
Since the form $i\t^1\w\ol\t^1\w  dz_2$ is harmonic, 
the class $[i\t^1\w\ol\t^1\w  dz_2]$ does not vanish.
\end{proof}

The surface $S_{t}^{(+)}$ admits a $1$-dimensional family of 
flat line bundles since $b_1(S_t^{(+)})$ is equal to $1$.
Tricerri's example of $l.c.k$ metric gives the corresponding flat line bundle $L$ which coincides with the canonical line bundle $K$. 
It is natural to ask whether there is an l.c.k metric which gives a different corresponding flat line bundle. However we have
\bgn{proposition}\label{prop:7.3}
Let $\ome'$ be an l.c.k structure on the Inoue surface $S:=S_{N,p,q,r;t}^{(+)}$. 
Then the corresponding flat line bundle $L'$ to $\ome'$ must be the canonical line bundle.
\end{proposition}
\bgn{proof}
Since the Inoue surface $S$ admits no curves, 
$H^0(S, K \otimes (L')^{-1})=\{0\}$.
If $\ome'$ gives a flat line bundle $L' \neq K$, 
then $H^2(S, L')\cong H^0(S, K \otimes (L')^{-1})=\{0\}$.
Then we can apply the stability theorem to obtain a family of l.c.k structures 
$\{\ome_t\}$, where $t\in \C$ is a parameter of deformations of complex structures. 
However it follows from Belgun's result that there is no such family of l.c.k forms. Thus $L'$ must be $K$. 
\end{proof}

\subsection{Inoue surfaces $S_M$ }
The Inoue surfaces $S_M$ are also the quotient surfaces $S_M=\Bbb H \times \C/ G_M$, 
where $\H$ is the upper half plane of the complex numbers $\C$ and $G_M$ is a group of 
analytic automorphisms of $\H\times\C$. 
Let $M \in$ SL$(3, \Z)$ be a unimodular matrix with a real eigenvalue $\a>1$ and two complex conjugate 
eigenvalues $\b \neq \ol\b$. 
Let $(a_1, a_2, a_3)$ be a real eigenvector corresponding to $\a$ and $(b_1, b_2, b_3)$ a eigenvector 
corresponding to $\b$. 
Then $G_M$ is the group generated by the following automorphisms: 
\bgn{align}
&\phi_0(w, z) \mapsto (\a w, \b z)\\
&\phi_i(w,z)\mapsto (w+a_i, z+b_i), \quad \text{\rm for } i=1,2,3.
\end{align}
Then the action of $G_M$ is properly discontinuous and has no fixed points. 
Thus $S_M=\H\times\C/G_M$ is a compact complex surface which is differentially a fibre bundle over the circle $S^1$ with the $3$-torus as fibre. 
On $S_M$, an l.c.k form $\ome$ is defined by 
$$
\ome=-i\(\frac{dw\w d\ol w}{(w_2)^2}+w_2dz\w d\ol z\).
$$
A basis of invariant forms of type $(1,0)$ is given by 
$$
\t^1=\frac{dw}{w_2}, \,\, \t^2=(w_2)^{\frac12}dz
$$
Then $\ome$ is written as 
$$
\ome=-i\(\t^1\w\ol\t^1+\t^2\w\ol\t^2\).
$$
By using  
$$
d\t^1=\frac1{2i}\ol\t^1\w\t^1,\qquad d\t^2=\frac12\frac{\t^1-\ol\t^1}{2i}\w\t^2,
$$
we have
\bgn{align}
d(\t^1\w\ol\t^1)=0, \quad d(\t^2\w\ol\t^2)=\eta\w\t^2\w\ol\t^2,
\end{align}
Thus $\ome$ is an l.c.k structure with Lee form $\eta=\frac{dw_2}{w_2}.$

We put  $d_\eta=d-\eta$. Then we have $d_\eta\ome=0$.
As in the previous section, we have the complex by using the operator $d_\eta$: 
$$
\cdots\overset{d_\eta}\arrow \w^k \overset{d_\eta}\arrow \w^{k+1}
\overset{d_\eta}\arrow\cdots
$$
The cohomology group of the complex is denoted by $H^k_\eta(S_M)$. Then we already see that
$H^k_\eta(S_M)\cong H^k(S_M, L)$.
Then we obtain 
\bgn{lemma}
$$[\eta\w\ome]\neq 0 \in H^3_{\eta}(S_M)\cong H^3(S_M,L).$$
\end{lemma}
\bgn{proof}
We apply the Hodge theory to the complex $(\w^\bullet, d_{\eta})$. 
Let $d_\eta^*$ be the formal adjoint operator of $d_\eta$ with respect to the Hermitian form $\ome$. 
Then $d_\eta^*$ is given by 
$$
d_\eta^*=(-1)^k*^{-1}(d_{-\eta})*=(-1)^k*^{-1}(d+\eta)*
$$
where $d^*_\eta$ acts on $k$-forms.
Note that $\eta$ changes into $-\eta$. 
Then we obtain 
\bgn{align}
*(\eta\w\ome)=&*(\frac{dw_2}{w_2}\w\ome)=-\frac{dw_1}{w_2}\\
d_{-\eta}*(\eta\w\ome)=&-(d+\eta)\frac{dw_1}{w_2}\\
=&\frac{dw_2}{w_2}\w\frac{dw_1}{w_2}-\frac{dw_2}{w_2}\w\frac{dw_1}{w_2}=0
\end{align}
Thus $\eta\w\ome$ is a harmonic form and the class $[\eta\w\ome] \in H^3_\eta(S_M)$
does not vanish.
\end{proof}
Thus we apply Theorem \ref{th:obstruction H3} and obtain
\bgn{proposition}\label{prop:7.5}
Let $\ome$ be the l.c.k structure on $S_M$ and $L$ the corresponding flat line bundle to $\ome$. 
We denote by $\{L_s\}$ deformations of flat line bundles with
$L_0=L$ which is given by $d$-closed $1$-forms $\{\eta s\}$, where $[\eta]\neq 0 \in H^1(S_M)$.
Then 
 $S_M$ does not admit a smooth family of l.c.k structures 
 $\{\ome_s\}$ such that $L_s$ is the corresponding flat line bundle to  $\ome_s$.
\end{proposition}

\subsection{Inoue surfaces $S_{N,p,q,r}^{(-)}$} 
As in the previous section, the surfaces $S_{N,p,q,r}^{(-)}$ are defined 
as quotient complex manifolds $\H \times \C/G_{N,p,q,r}^{(-)}$. 
The group $G_{N,p,q,r}^{(-)}$ is generated by the following automorphisms:
\bgn{align}
&\phi_0:(w,z) \mapsto (\a w, -z), \\
&\phi_i:(w,z) \mapsto (w+a_i, z+b_i w+ c_i)\\
&\phi_3:(w,z)\mapsto \(w, z+\frac{(b_1a_2-b_2a_1)}{r}\),
\end{align}
where real constants $\a, a_i, b_i, c_i, r$ are the same as in the subsection of $S_{N,p,q,r;t}^{(+)}$.
Then we have a following basis of forms of type $(1,0)$ on $\H\times\C$:
$$
\t^1=\frac{dw}{w_2}, \quad \t^2=\frac{z_2}{w_2}dw-dz
$$
The forms $\t^1$ and $\t^2$ are invariant under the action of $\phi_i$ and 
$\phi_3$ for $i=1,2$ and $\phi^*_1\t^1=\t^1$ and $\phi^*_1\t^2=-\t^2$. 
Thus $\ome=-i(\t\w\ol\t^1+\t^2\w\ol\t^2)$ is an invariant form which is an l.c.k structure on $S_{N,p,q,r}^{(-)}$, that is,
$$
d \ome=\frac{dw_2}{w_2}\w\ome.
$$
The l.c.k form $\ome $ gives the corresponding flat line bundle $L$. Deformations of flat line bundles $\{ L_s\}$ with $L_0=L$ are given by a class $[s \eta]\in H^1(S_{N,p,q,r}^{(-)}, \R)$. 
Then we also have 
\bgn{proposition}
The inoue surface $S_{N,p,q,r}^{(-)}$ does not admit a smooth family of l.c.k structures
 $\{\ome_s\}$ such that $L_s$ is the corresponding flat line bundle to  $\ome_s$.
\end{proposition}\label{prop:7.6}
\bgn{proof}
We denote by $S^{(-)}$ the Inoue surface $S^{(-)}_{N,p,q,r}$.
It suffices to show that the class $[\eta\w\ome]\in H^3_\eta(S^{(-)})$ does not vanish. 
We see that 
$d_\eta(\eta\w\ome)=0$. The formal adjoint $(d_{\eta})^*$ is given by 
$$
(d^*_{\eta})^*=(-1)^k*^{-1}\circ d_{-\eta}\circ *=(-1)^k*^{-1}\circ (d+\eta)\circ *,
$$
where $(d_\eta)^*$ acts on $k$-forms and $\eta=\frac{dw_2}{w_2}$.
Since $*(\eta\w\ome)=-\frac{dw_1}{w_2}$, we have 
\bgn{align}
d_{-\eta}\circ *(\eta\w\ome)=&-d_{-\eta}\frac{dw_1}{w_2}\\
=&\frac{dw_2}{w_2}\w\frac{dw_1}{dw_2}-\frac{dw_2}{w_2}\w\frac{dw_1}{dw_2}=0
\end{align}
Thus $d_{\eta}^*(\eta\w\ome)=0$.
Thus $\eta\w\ome$ is harmonic and the class $[\eta\w\ome]\in H^3_\eta(S^{(-)})$
does not vanish.
\end{proof}

\bgn{thebibliography}{99}
\bibitem{Bel}
F.~A.~Belgun, {\it On the metric structure of non-K\"ahler complex surfaces}, Math. Ann. 317 (2000), no. 1, 1--40, MR1760667.

\bibitem{BG}
P.~C.~ Boyer and K.~ Galicki,
{\it Sasakian geometry},
Oxford Mathematical Monographs; Oxford Science Publications. Oxford: Oxford University Press (2008).

\bibitem{Bru1}
 M.~Brunella, {\it Locally conformally K\"ahler metrics on certain non-K\"ahlerian surfaces}, Math. Ann. 346 (2010), no. 3, 629--639, MR2578564.
\bibitem{Bru2}
M.~Brunella, {\it Locally conformally K\"ahler metrics on Kato surfaces},
arXiv:1001.0530 
\bibitem{DO}
 S.~Dragomir and L.~Ornea, {\it Locally conformal K\"ahler geometry. Progress in Mathematics}, 155. Birkh\"auser Boston, Inc., Boston, MA, 1998. xiv+327 pp. ISBN: 0-8176-4020-7,
 MR1481969 (99a:53081).
\bibitem{FP}
A. Fujiki and M. Pontecorvo, 
{\it Anti-self-dual bihermitian structures on Inoue surfaces}, 
J. Differential Geometry, 85 (2010), no. 1, 15--72,

\bibitem{Go-1}
R.~Goto,
{\it Moduli spaces of topological calibrations,
Calabi-Yau, hyperK\"ahler, G$_2$ and Spin$(7)$ structures},
Internat. J. Math. 15 (2004), no. 3, 211--257. MR 2060789, Zbl 1046.58002
\bibitem{Go0}
R.~Goto,
{\it On deformations of generalized Calabi-Yau, hyperK\"ahler, G$_2$ and Spin$(7)$ structures},
Math.DG/0512211
\bibitem{Go1}
R.~Goto,
{\it Deformations of \complex  and generalized K\"ahler structures},
 J. Differential Geometry, 84 (2010), no. 3, 525--560,
Math. DG/0705.2495
\bibitem{Go2}
R.~Goto,
{\it Poisson structures and generalized K\"ahler submanifolds},
 J. Math. Soc. Japan 61 (2009), no. 1, 107--132. MR 2272873, Zbl 1160.53014
 \bibitem{Go3}
 R.~Goto,
 {\it Deformations of Generalized K\"ahler structures and Bihermitain structures},
 Math. DG/0910.1651
 \bibitem{GO}
 P.~Gauduchon and L.~ Ornea, {\it Locally conformally K\"ahler metrics on Hopf surfaces},
 Ann. Inst. Fourier (Grenoble) 48 (1998), no. 4, 1107--1127, MR1656010. 
 \bibitem{Ino}
 M.~Inoue, {\it On surfaces of Class ${\rm VII}_{0}$},
  Invent. Math. 24 (1974), 269--310, MR0342734.
  \bibitem{Ise}
M.~Ise, {\it On the geometry of Hopf manifolds}, Osaka Mathematical Journal. Volume 12, Number 2 (1960), 387-402. 
\bibitem{Ko}
K.~Kodaira,
{\it Complex manifolds and deformations of complex structures},
Grundlehren der Mathematischen Wissenschaften [Fundamental Principles of Mathematical Sciences], 283. Springer-Verlag, New York, 1986. x+465 pp. ISBN: 0-387-96188-7 MR 0815922, Zbl 0581.32012
Grundlehren der Mathematischen Wissenschaften, {\bf 283}, springer-Verlag,
(1986)
\bibitem{K.S.I,II}
K.~Kodaira and D.C.~Spencer,
{\it  On deformations of complex, analytic structures I,II},
Ann. of Math. (2) 67 (1958) 328--466. MR 0112154, Zbl 0128.16901
\bibitem{K.S III}
K.~Kodaira and D.C.~Spencer,
{\it On deformations of complex analytic structure, III.
stability theorems for complex structures},
Ann. of Math. (2) 71 (1960) 43--76. MR 0115189, Zbl 0128.16902 
\bibitem{LLMP}
M.~de Le\'on, B. L\'opez, J.C.~Marrero and E.~Padr\'on, 
{\it On the computation of the Lichnerowicz-Jacobi cohomology}, 
J. Geom. Phys. 44 (2003), no. 4, 507--522, MR1943175.
\bibitem{Mall}
D.~ Mall, {\it The cohomology of line bundles on Hopf manifolds}, Osaka Journal of Mathematics. Volume 28, Number 4 (1991), 999-1015. 

\bibitem{OV1}
L.~Ornea and M.~ Verbitsky,
{\it Locally conformal K\"ahler manifolds with potential},
 Math. Ann. 348 (2010), no. 1, 
MR2657432.
\bibitem{OV2}
 L.~Ornea and M.~ Verbitsky, {\it Morse-Novikov cohomology of locally conformally K\"ahler manifolds}, J. Geom. Phys. 59 (2009), no. 3, 295--305, MR2501742.
 \bibitem{Tri} 
F.~Tricerri, 
{\it Some examples of locally conformal K\"ahler manifolds},
Rend. Sem. Mat. Univ. Politec. Torino 40 (1982), no. 1, 81--92,
MR0706055. 
\bibitem{Vai}
I.~Vaisman, {\it Generalized Hopf manifolds}, Geom. Dedicata 13 (1982), no. 3, 231--255, MR069067.
\end{thebibliography}
\end{document}